\documentclass[11pt]{article}
\usepackage{fleqn}
\usepackage{array}
\usepackage[latin1]{inputenc}
\usepackage{amssymb}
\usepackage{amsfonts}
\advance \topmargin by -\headheight
\advance \topmargin by -\headsep
\evensidemargin \oddsidemargin
\newcommand{\bE}{{\mathbb{E}}}
\newcommand{\D}{{\displaystyle}}

\newcommand{\cC}{{\cal C}}
\newcommand{\cG}{{\cal G}}
\newcommand{\cI}{{\cal I}}
\newcommand{\cF}{{\cal F}}
\newcommand{\cO}{{\cal O}}
\newcommand{\cM}{{\cal M}}
\newcommand{\cJ}{{\cal J}}
{

\setlength{\parskip}{10pt plus 2pt minus 1pt}
\setlength{\parindent}{0pt}

\setcounter{section}{-1}

\begin{document}
\begin{center}
\section*{Coincidences of simplex centers and related facial structures}
\end{center}

\begin{center}
Allan L. Edmonds$^a$, Mowaffaq Hajja$^b$, Horst Martini$^c$\footnote{The third
author was partially supported by the Yarmouk University in Irbid (Jordan) and by
a DFG grant.}

$^a$ Department of Mathematics\\Indiana
University\\Bloomington, IN 47405\\USA\\[0.3cm]
$^b$ Department of Mathematics\\Yarmouk
University\\Irbid\\JORDAN\\[0.3cm]
$^c$ Faculty of Mathematics\\Chemnitz University of
Technology\\09107 Chemnitz\\GERMANY
\end{center}

\begin{abstract}
We investigate geometric properties of simplices in the Euclidean
$d$-dimensional space for which analogues of the classical
triangle centers coincide. We try to give a unified presentation
of related results, partially also unifying known results
for $d=2$ and $d=3$.\\[0.5cm]

\textbf{Keywords:} barycentric coordinates, cevian, Cholesky decomposition, Gram matrix, (regular) simplex,
simplex center
\end{abstract}

\section{Introduction}

Mainly by methods from linear algebra, we study the analogues of the classical triangle centers (cf. \cite{Ki})
for general simplices in the Euclidean $d$-dimensional space $\bE^d, d
\ge 2$. We focus on interpreting the significance of two or more
of the classical centers coinciding. Also we give several
instructive constructions of examples. The centers under study
include the centroid, the circumcenter, the incenter, the
orthocenter (or its proper higher dimensional generalization, the
Monge point), and  the Fermat-Torricelli point. We also consider two
new centers with clear geometric meanings. In the last section we examine classes of
simplices whose cevians through certain centers have equal lengths.

In dimension $d=2$ it is a standard fact that if two classical
centers coincide, then the triangle is equilateral, see, e.g., \cite{Isa}.

In dimension $d=3$, the parallel conclusion is that if two
classical centers coincide, then the tetrahedron is equiareal
(i.e., it has faces of equal area implying that these
faces are even congruent), but not necessarily equilateral
(regular). One noteworthy point is that the orthocenter does not
necessarily exist for $d \ge 3$. But when it does, one can often
make much stronger conclusions, both for $d=3$ and in higher
dimensions. We hope to return to this point in a subsequent paper,
where we plan to make a detailed analysis of orthocentric
simplices.

In dimension $d \ge 4$ the situation becomes yet more complicated.
When two of the classical centers coincide one can give a
meaningful geometric description in terms of the facial structure of the
studied simplices. But various examples of
various degrees of subtlety show that in general, when two centers
coincide, one cannot usually infer much about other centers. After
having done some of the work described here we learned of the
not-well-known papers of V. Devide (see \cite{De 1} and \cite{De 2}) where he also treated
some of these problems. We give our own proofs, however. We also
resolve questions he posed but did not answer.

Some of the material discussed here exists in various forms in
older, scattered literature. We intend to give a unified
presentation, collecting a number of related results, with proofs
as well as references to the known literature. We also would like to mention the papers
\cite{Ed}, \cite{F-M}, \cite{McM}, \cite{We} and \cite{M-W}, where special types of simplices
are investigated, but only regarding their facial structure and not in view of
coincidence of certain centers; see also the survey \cite{Ma}, \S~9. And in \cite{Ave} some
related results are given for simplices in normed linear spaces.

\section{Terminology and notation}

A \emph{$d$-simplex} $S = [A_1, \dots, A_{d+1}]$ in the Euclidean
$d$-dimensional space $\bE^d, \, d \ge 2$, with origin ${\bf o}$ is
defined as the convex hull of $d+1$ affinely independent points
(or position vectors) $A_1, \dots, A_{d+1}$ in $\bE^d$. Thus the
vectors $A_i - A_j, \, 1 \le i \le d+1, \, i \not= j$, are
linearly independent for every $j \in \{1, \dots, d+1\}$, and
therefore the linear dependence relation
\[
c_1 A_1 + \dots c_{d+1} A_{d+1} = {\bf 0}
\]
is unique up to multiplication by a scalar.
The points $A_1, \dots, A_{d+1}$ are called the \emph{vertices} of
$S$. A line segment that joins two vertices of $S$ is called an
\emph{edge} of $S$, and a $j$-simplex whose vertices are any $j+1$
vertices of $S$ is said to be a \emph{$j$-face} of $S$. The
$(d-1)$-simplex whose vertices are all vertices of $S$ except for
$A_j, \, j \in \{1, \dots, d+1\}$, is called the $j$th
\emph{facet} of $S$, or the \emph{facet opposite to} $A_j$. The
$d$-simplex $S$ is \emph{regular} if all its edges have
equal length, it is \emph{equiareal} if all its facets
have the same $(d-1)$-volume, and it is called \emph{equifacetal}
if all its facets are congruent (i.e., isometric, see \cite{Ed} for
interesting new results on equifacetal simplices). Moreover, we
say that $S$ is \emph{equiradial} if all its facets have the same
circumradius (see (4) below for a definition of the circumradius).

Every point $P$ in the convex hull of $A_1, \dots, A_{d+1}$ can
be represented in the form
\begin{equation}
P = \frac{v_1 A_1 + \dots + v_{d+1} A_{d+1}}{v}\,,
\end{equation}
where $v_j$ is the $d$-volume of the $d$-simplex obtained from $S$ by
replacing $A_j$ by $P$. So the $d$-volume $v$ of $S$ is represented by
$v = v_1 + \dots + v_{d+1}$. Note also that $v_j = \frac{1}{d} a_j
h_j$, where $a_j$ is the $(d-1)$-volume of the $j$th facet of $S$,
and $h_j$ is the altitude from  $A_j$ to the $j$th facet, see,
e.g., \cite{Be}, Theorem 9.12.4.4, page 260. The numbers
$\frac{v_j}{v}$ in (1) are called the \emph{barycentric
coordinates} of $P$ with respect to $A_1, \dots, A_{d+1}$.

The \emph{centroid} $\cG$ of the $d$-simplex $S = [A_1, \dots,
A_{d+1}]$ is defined as the average
\begin{equation}
\cG = \frac{A_1 + \dots + A_{d+1}}{d+1}
\end{equation}
of its vertices.

The \emph{insphere} of $S$ is the sphere that is tangent to all
$d+1$ facets of $S$; its center is the \emph{incenter} $\cI$ of $S$,
and its radius is the \emph{inradius} of $S$. Since $\cI$ is
equidistant to all the facets of $S$, it follows that the $j$th
barycentric coordinate of $\cI$ is proportional to the
$(d-1)$-volume of the $j$th facet, i.e., the incenter is algebraically
defined by
\begin{equation}
\cI = \frac{a_1 A_1 + \dots + a_{d+1} A_{d+1}}{a_1 + \dots +
a_{d+1}}\,,
\end{equation}
where $a_i$ is the volume of the $i$th facet of $S$.
The \emph{circumsphere} of $S$ is the sphere passing through all
vertices of $S$, and the center $\cC$ of that sphere is called the
\emph{circumcenter} of $S$. Thus $\cC$ is defined by the
requirement that
\begin{equation}
|| \cC - A_i || = || \cC - A_j || \, \mbox{ for } \, 1 \le i \le j \le
d+1\,,
\end{equation}
where $||\cC-A_i||$ is said to be the \emph{circumradius} of $S$.
Unlike the centroid and the incenter, the circumcenter may lie
outside of $S$.

The \emph{Fermat-Torricelli point} $\cF$ of $S$ is defined to be the
point whose distances to all the vertices of $S$ have minimal sum.
Such a point exists and is unique, and setting
\[
f(i) = \sum \left( \frac{A_j - A_i}{||A_j - A_i||}: 1 \le j \le
d+1, \, j \not= i\right)||
\]
it is known (see \cite{KM 2}, Theorem 1.1, and \cite{BMS}, Theorem
18.3 and Reformulation 18.4, cf. further also \cite{Da}) that if $||f(i)|| > 1$ for all $i \in
\{1, \dots, d+1\}$, then
$\cF$ is an interior point of $S$ (floating case), and that if $||f(i)|| \le 1$ for some
$i$ then this $i$ is unique and $\cF = A_i$ (absorbed case). In the floating case,
$\cF$ is characterized by the property
\begin{equation}
\sum \limits^{d+1}_{i=1} \frac{\cF-A_i}{||\cF-A_i||} = {\bf 0}\,.
\end{equation}
For an interesting application of the Fermat-Torricelli point in classical
geometry, namely an extension of Napoleon's theorem to $d$-dimensional space,
we refer to \cite{Ma-Wei}.

The \emph{orthocenter} $\cO$ of a
$d$-simplex $S$ is, if it exists, the intersection of the $d+1$
altitudes of $S$. In stark contrast to the case $d=2$, a
$d$-simplex might not have an orthocenter when $d \ge 3$, see
\cite{AC} for $d=3$ and \cite{H-W}, \cite{Man} as well as
\cite{Mo} for higher dimensions.

The last center we want to define here is discussed only in the next section, i.e., for
$d=3$ (but we introduce it for $d$ arbitrary). For each edge $E_{ij} = A_i A_j$ of $S = [A_1,
\dots, A_{d+1}]$ there is a unique hyperplane $H_{ij}$ containing the centroid $\cG_{ij}$
of the remaining $d-1$ vertices and perpendicular to $E_{ij}$. These ${d+1 \choose 2}$ hyperplanes
have a common point, the \emph{Monge point} $\cM$ of $S$. This point $\cM$ is a reflection of
$\cC$ in $\cG$ and coincides, if $S$ is orthocentric, with the orthocenter $\cO$, see
also \cite{Mo}.

\section{Tetrahedra whose centers coincide}

Any $3$-simplex (or non-degenerate tetrahedron) has the centers
$\cG, \cI, \cC, \cF$, and $\cM$ (for $\cF$ see \cite{KM 1}, and for $\cM$ see
\cite{AC}, Article 229, pp. 76-77). It follows immediately that
the Monge point of the tetrahedron $S = [A,B,C,D]$ coincides with its orthocenter if $S$ is orthocentric,
and that this is equivalent to

\begin{equation}
A \cdot B = C \cdot D\,, \,\, A \cdot C = B \cdot D\,, \,\, A \cdot D
= B \cdot C\,,
\end{equation}
where ``$\cdot$'' means the ordinary inner product. Only parts of the following theorem can be found
in the basic references \cite{AC}, \cite{CB}, \cite{Th}, and \cite{Za}, which collect geometric
properties of tetrahedra in the Euclidean $3$-space.

\textbf{Theorem 2.1:} \emph{For a tetrahedron $T \subset \bE^3$
the following conditions are equivalent}:
\begin{enumerate}
\item \emph{The tetrahedron $T$ is equifacetal}.
\item \emph{The tetrahedron $T$ is equiareal}.
\item \emph{Every two opposite edges of ~$T$ are equal}.
\item \emph{The perimeters of the facets of ~$T$ are equal}.
\item \emph{The circumradii of the facets of~ $T$ are equal, i.e., $S$ is equiradial}.
\item \emph{The centroid, the incenter, the circumcenter, the Fermat-Torricelli
point and the Monge point of ~$T$ coincide}.
\item \emph{Two of the five centers mentioned above coincide}.
\end{enumerate}

\textbf{Proof:} It is clear that 1. implies all the other
statements, because the group of isometries of an equifacetal tetrahedron $[A,B,C,D]$
is transitive, as it contains the Klein $4$-group consisting of the permutations $\{(A\,B)\,(C\,D),(A\,C)\,
(B\,D), (A\,D)\,(B\,C),e\}$
(see \cite{Ed}, where this is proved even for any dimension).
That 2. implies 1. is the well-known theorem usually
referred to as Bang's Theorem (cf. \cite{Ha 1},
\cite{Ho}, pp. 90-97, \cite{Br}, \cite{AC}, Article 306, p.
108). The implication 3. $\Rightarrow$ 1. is trivial, 4. implies
3. by solving the corresponding system of linear equations, and 5.
implies 1. by \cite{Ha-Wa 1}, Theorem 3. That 6. implies 1. follows
from \cite{Ha-Wa 1}, Theorem 5, which states that if any two of $\cG,
\cI, \cC$, and $\cF$ coincide, then the tetrahedron $T$ is equifacetal.
This also follows from \cite{AC}, Article 305, page 108, which
states that if any two of $\cG, \cI, \cC$, and $\cM$ coincide, then $T$ is
equifacetal, cf. also \cite{De 2}. Thus it remains to show that 7.
implies 6., and in view of \cite{Ha-Wa 1}, Theorem 5, and \cite{AC},
Article 305, page 108, it suffices to show that if $\cM$ and $\cF$
coincide, then $T$ is equifacetal. So assume $\cM = \cF = {\bf 0}$, and that
${\bf 0}$ is in the interior of $T = [A,B,C,D]$. Set
\[
||A|| = \frac{1}{a}\,,\, \, ||B|| = \frac{1}{b}\,, \,\, ||C|| =
\frac{1}{c}\,, \,\, ||D|| = \frac{1}{d}\,.
\]
Then $aA + bB + cC + dD = {\bf 0}$, yielding by multiplication with $aA$
that if $A \cdot B = A \cdot C = 0$, then $(aA) \cdot (dD) = -1$, and
so ${\bf 0}$ would be on the line $AD$, contradicting the assumption
that $o$ is an interior point of $T$. Hence at most one of $A
\cdot B, A \cdot C$, and $A \cdot D$ is zero, and we may assume
that $A \cdot B \not= 0$ and $A \cdot C \not= 0$. Taking norms of
both sides of
\[
aA + bB = - cC - dD\,,
\]
we obtain $abA \cdot B = cd C \cdot D$. Since $0 \not= A \cdot B =
C \cdot D$, it follows that $ab = cd$, and similarly $ac = bd$.
Hence $ab = cd = ac = bd$, each being equal to $\sqrt{abcd}$.
Therefore $a=b = c = d$, and ${\bf 0}$ is the centroid $\cG$ of $T$. Thus
$\cG$ coincides with the Fermat-Torricelli point of $T$, and by (5)
$T$ is equifacetal. \hfill $\Box$

\textbf{Remark 1:} It should be mentioned that equifacetal tetrahedra can be used to
give interesting characterizations of Euclidean motions, see
\cite{Le}. Also we mention here that equifacetal/equiareal tetrahedra are called \emph{isosceles tetrahedra}
by many authors, see, e.g., \cite{AC}. We will not follow that way, since we use the notion
of \emph{isosceles simplices} in another sense, see the proofs of Theorems 3.3, 3.4, and in Section 4 below.

Motivated by 5. in Theorem 2.1 one might ask whether equality of
the inradii of the facets of a tetrahedron implies equifacetality.
Also it should be interesting to check whether all remains valid if
more centers are added to the list in 6. of Theorem 2.1. Negative
answers to both these questions are supplied in Theorems 2.2 and
2.3 below.

\textbf{Theorem 2.2:} \emph{The inradii of the facets of the
non-equifacetal tetrahedron, whose edges have lengths $1,1,1,1,1,
(3 + \sqrt{33})/6$, are equal}.

\textbf{Proof:} To justify the existence of a tetrahedron whose
edges are as given, and whose edge lengths are, more generally, even
equal to $1,1,1,1,1, t$ with $t \in (0, \sqrt{3})$, we start with
a rhombus $ABCD$ whose sides all have unit length and whose short
diagonal is $AC$. Keeping $ABC$ fixed, we fold $ABCD$ against
$AC$, letting $D$ move towards $B$. The tetrahedra $T_t
= [A,B,C,D]_t$ formed in this way have edges of length
$1,1,1,1,1,t$ with $t$ ranging in $(0,\sqrt{3})$. Now let $r=r(t)$
be the inradius of a triangle whose side-lengths are $1,1,t$.
Then $r= \frac{2a}{p}$, where $a$ is the area and $p$ the
perimeter of the triangle. Heron's formula (cf. \cite{Co}, \S~1.5) yields
\[
f(t) : = 4 (r(t))^2 = \frac{16a^2}{p^2} = \frac{t^2(2-t)}{t+2}\,,
\]
and solving $f(t) = f(1)$ we find that $t = (3+\sqrt{33})/6$, as
desired. \hfill $\Box$

It is also interesting to find new natural centers whose coincidence with known ones does not imply
that a tetrahedron is equifacetal. For this purpose we set
\[
\cJ = \frac{\alpha A + \beta B + \gamma C + \delta D}{p}
\]
for a tetrahedron $T = [A,B,C,D]$, where $\alpha, \beta, \gamma, \delta$ are the perimeters of the
facets of $T$ opposite to $A,B,C,D$, respectively, and $p = \alpha
+ \beta + \gamma + \delta$. We call $\cJ$ the \emph{complementary
$1$-centroid} of $T$.

\textbf{Theorem 2.3:} \emph{The incenter and the complementary
$1$-centroid of a tetrahedron coincide iff the inradii of its
facets are equal. Consequently, there exist
non-equifacetal tetrahedra $T$ whose incenter and complementary
$1$-centroid coincide}.

\textbf{Proof:} By (3) (replacing $a_1, \dots, a_4$, i.e., the areas of the
faces of $T$, by $a,b,c,d$)
and the definition of the complementary $1$-centroid the points
$\cI$ and $\cJ$ coincide iff
\[
\frac{a}{\alpha} = \frac{b}{\beta} = \frac{c}{\gamma} =
\frac{d}{\delta}\,\,,
\]
which is equivalent to the property that the inradii of the facets
are equal, since $\frac{2a}{\alpha}$ is the inradius of the facet
opposite to $A$, etc.; the latter statement follows from Theorem
2.2. \hfill $\Box$

\section{Higher-dimensional simplices whose centers coincide}

The situation in $3$-space does not have exact analogues in higher
dimensions. In the following we exhibit various related results
for $d$-simplices if $d \ge 4$ or, in some cases, if $d \ge 3$.

\textbf{Theorem 3.1:} \emph{Let $S = [A_1, \dots, A_{d+1}]$ be a
$d$-simplex. If any two of the centroid, the circumcenter
and the Fermat-Torricelli point of $S$ coincide, then all
three centers coincide}.

\textbf{Proof:} From the definitions it follows that
\[
\begin{array}{lll}
{\bf 0} = \cG  & \Leftrightarrow & A_1 + \dots + A_{d+1} = {\bf 0}\,,\\
{\bf 0} = \cF & \Leftrightarrow & \frac{A_1}{||A_1||} + \dots +
\frac{A_{d+1}}{||A_{d+1}||} = {\bf 0}\,,\\
{\bf 0} = \cC & \Leftrightarrow & ||A_1|| = \dots = ||A_{d+1}||\,.
\end{array}
\]
Since the dependence relation among $A_1, \dots, A_{d+1}$ is
unique up to multiplying by a scalar, the proof is complete. \hfill
$\Box$

We remark that the four statements of the following theorem were proven in \cite{De 1}.
However, we give partially new, shorter proofs. And also we need an additional notion, yet. Namely,
we say that a $d$-simplex $S$ has \emph{well-distributed edge lengths} if \emph{all} its
facets have the same sum of squares of
all their ${d \choose 2}$ edge lengths.

\textbf{Theorem 3.2:} \emph{For any $d$-simplex $S$ the following statements hold true}.
\begin{enumerate}
\item[(i)] \emph{The centroid $\cG$ and the circumcenter $\cC$ of $S$ coincide iff $S$
has well-distributed edge lengths}.
\item[(ii)] \emph{The circumcenter $\cC$ and the incenter $\cI$ of $S$ coincide iff $S$ is
equiradial}.
\item[(iii)] \emph{The incenter $\cI$ and the centroid $\cG$ of $S$ coincide iff $S$ is
equiareal}.
\item[(iv)] \emph{The points $\cG,\cC$, and $\cI$ of $S$ coincide iff two of the three conditions
$\{S$ has well-distributed edge lengths;  $S$ is equiradial;  $S$ is
equiareal}$\}$ \emph{hold, in each case implying the third}.
\end{enumerate}

\textbf{Proof:} For the proofs of statements (i) and (iv) we refer to \cite{De 1}, the
other two equivalences will now be verified by proofs which are shorter than those given in \cite{De 1}.
To see (ii), drop a perpendicular from $\cC$ to the facet opposite to $A_i, \, i \in
\{1, \dots, d+1\}$. The obtained
intersection point is the circumcenter $\cC_i$ of the $i$th facet. The distance of $\cC$ to any
vertex of $S$ is the circumradius $R$ of $S$, and the distance from $\cC_i$ to any vertex of $S$
different from $A_i$ is the circumradius $R_i$ of the $i$th facet. So the three points $\cC, \cC_i,
A_j \, (j = 1, \dots, d+1; i \not= j)$ form a right triangle, and $R^2 = R^2_i + |\cC\cC_i|^2$.
But if $\cI = \cC$, then $|\cC\cC_i|^2 = r^2$ (the squared inradius of $S$). Hence $R^2_i = R^2 - r^2$,
not depending on the choice of the facet, and $S$ is equiradial. On the other hand, if $R^2_i$ does
not depend on the choice of the facet, then the formula yields $|\cC\cC_i|^2 = R^2 - R^2_i$,
also independent of the choice of the facet. Thus $\cC$ has to be the incenter of $S$.
Finally we show (iii). By (2) and (3) we see that ${\bf 0}$ is the centroid of $S$ iff $A_1 + \dots
+ A_{d+1} = {\bf 0}$, and that ${\bf 0}$ is the incenter of $S$ iff $v_1 A_1 + \dots + v_{d+1} A_{d+1} = {\bf 0}$,
where $v_i$ is the $(d-1)$-volume of the $i$th facet of $S$ opposite to $A_i$. Since the
dependence relation among $A_1, \dots, A_{d+1}$ is unique up to multiplying by a scalar, (iii)
is obtained. \hfill $\Box$

M. Hajja and P. Walker \cite{Ha-Wa 1} have shown that a tetrahedron satisfying
$\cF=\cI$ must be equifacetal, see also Theorem 2.1. It follows that an orthocentric
tetrahedron (i.e., a tetrahedron with orthocenter) in which $\cF=\cI$ must be regular. A similar statement
holds in dimension 2.
One may conjecture that in any dimension an orthocentric simplex
with $\cF=\cI$ must be regular. This certainly is consistent with
dimensional analysis: The space of orthocentric $d$-simplices is
$(d+1)$-dimensional. The equality $\cF=\cI$ amounts to $d$
equations, one for each coordinate. That leaves one degree of
freedom, which can be accounted for by scaling. Since now we will
show that a $d$-simplex, $d \ge 4$, with
$\cF=\cI$ need not equifacetal, any proof of the conjecture must
somehow mix the two hypotheses together more than is done in
dimension 3.

\textbf{Theorem 3.3:} \emph{For any $d \ge 4$, a $d$-simplex whose incenter $\cI$ and Fermat-Torricelli point
$\cF$ coincide need not be equifacetal}.

\textbf{Proof.} We say that a $d$-simplex is \emph{isosceles} if it has a vertex $P$ such that
all edges emanating from $P$ have the same length. (Note that this definition generalizes the
standard notion in dimension $2$, but differs from the use of the term for $d=3$ in some other sources,
such as in \cite{AC}; see Remark 1 in Section 2 above.) So we denote an isosceles $d$-simplex $S$ with base
$T$ and opposite vertex $P$ by $S = [T,P]$.
We will assume that $T$ is equifacetal. It is known that
non-regular equifacetal simplices exist in abundance arbitrarily
near any regular simplex. As an equifacetal simplex, $T$ has a unique
center, which we arrange to lie at the origin ${\bf 0} \in {\bf
\bE}^{d-1}$, where $T \subset {\bE}^{d-1} \subset {\bE}^d$. We
also assume that $P = ({\bf 0}, h)$.

Let $R$ and $r$ denote the circumradius and inradius of $T$. It is
known that $R \ge (n-1) r$, with equality if and only if $T$ is
regular.
Suppose $T = [A_1, \dots, A_d]$. Since $T$ is equifacetal, we know
that $\sum^d_{i=1} A_i = {\bf 0}$, since the centroid is ${\bf 0}$, and
also that $|A_1| = \dots = |A_d|$, since the circumcenter is ${\bf 0}$.
By symmetry all centers of $S$, such as the incenter and the
Fermat-Torricelli point, have the form $({\bf 0}, z)$. Let $\cF = ({\bf 0}, f)$
and $\cI = ({\bf 0}, i)$.
Now, provided $S$ is not too short, the Fermat-Torricelli point $\cF$ of $S$ is
characterized by the condition (5)
(where we note that $\frac{P-\cF}{||P-\cF||} = ({\bf 0},1))$.

Now $||A_i - \cF|| = \sqrt{R^2+f^2}$, so this sum yields
\[
d({\bf 0},-f) / \sqrt{R^2 + f^2} = ({\bf 0},-1) \, \mbox{ or }\, df = \sqrt{R^2 + f^2}\,.
\]
Hence
\[
(d^2 - 1) f^2 = R^2 \, \mbox{ or } \, f = \frac{R}{\sqrt{d^2-1}}\,.
\]
In particular $f$ and hence $\cF$ do not depend on $h$, provided $h$
is big enough, at least. (Otherwise $\cF=P$.) The
condition we need is just that
\[
h > \frac{R}{\sqrt{d^2-1}}\,.
\]
Now as $h$ increases from $0$ toward $\infty, i$ increases from
$0$ toward $r$. So we can choose $h$ so that $i =
\frac{R}{\sqrt{d^2-1}}$ provided that $\frac{R}{\sqrt{d^2-1}} <
r$. But we know that $R \ge (d-1)r$, with equality iff $T$ is
regular. So the possible region of success is
\[
(d-1)r < R < \left(\sqrt{d^2-1}\right) r\,.
\]
When $h = \frac{R}{\sqrt{d^2-1}}$, then $i$ is certainly much less
than $\frac{R}{\sqrt{d^2-1}}$, since the incenter lies in the
interior of $S$. As $h$ increases, $i$ also increases toward $r$.
If we choose $T$ to be equifacetal but not regular, very near the
regular $(d-1)$-simplex (as in \cite{Ed}, using $d \ge 4$), then $R$
is close to, but larger than $(d-1)r$, hence less than
$(\sqrt{d^2-1})r$.
For this purpose, it is convenient to rescale each simplex under
consideration so that $r=1$ for all $(d-1)$-simplices considered. \hfill $\Box$

\textbf{Theorem 3.4:} \emph{For any $d \ge 4$, there are equiradial $d$-simplices which are
not equiareal}.

\textbf{Proof:}  We first show that if $S = [T,P]$ is an isosceles $d$-simplex
with vertex $P$, base $T$ (i.e., facet opposite to $P$), and edge length $h$ at $P$, then the
circumradius $R_S$ of $S$ is given by
\[
R_S = \frac{h^2}{2\sqrt{h^2-R^2_T}}\,.
\]
To see this, let $\cC_S$ and $\cC_T$ denote the circumcenters of $S$ and $T$, respectively, and note that
$\cC_S$ lies on the line $P\cC_T$, which is perpendicular to $T$. Let $k$ denote the distance between
$\cC_S$ and $\cC_T$. Let $V$ be a vertex of $T$. Now, applying the Pythagorean theorem to the triangles
$P\cC_{T}V$ and $\cC_{S}\cC_{T}V$, with right angles at $\cC_T$, we have $(R_{S} \pm k)^2 + R_{T}^2 = h^2$ and
$R_{S}^2 - R_{T}^2 = k^2$. (Use a plus sign if $\cC_S$ lies between $P$ and $\cC_T$, and a minus sign if $\cC_T$
lies between $P$ and $\cC_S$.) One may solve the second equation for $k$ and substitute the result
in the first equation, getting
\[
\left(R_{S} \pm \sqrt{R_{S}^2-R_{T}^2}\right)^2 + R_{T}^2 = h^2\,.
\]
Expanding and collecting terms yields $\pm 2 R_{S} \sqrt{R_S^2-R_T^2} = h^2 - 2R_S^2$. Squaring both
sides, we have
\[
4R_S^4(R_S^2-R_T^2) = h^4-4h^2R_S^2+4R_S^4\,,
\]
hence $4(h^2-R_T^2) R_S^2 = h^4$, which yields the desired formula for $R_S$.
Now we show that for any $d \ge 4$ there is an isosceles
$d$-simplex with an equilateral (or a regular) base that is equiradial but not equiareal. Let $T$ be a regular
$(d-1)$-simplex, normalized for convenience to have edge length $1$, say. For isosceles
$d$-simplices of the form $[T,P]$, where $P$ has distance $h$ to each of the vertices of $T$,
we need to calculate what values of $h$ yield equiradial $d$-simplices. Certainly, $h = 1$ works
in any dimension, producing the regular $d$-simplex. But when $d \ge 4$, there is a second
value of $h$ yielding the desired examples.
Note that all facets $F$ of $T$ have the same circumradius. We seek an edge length $h$
such that $R_{[F,P]} = R_T$ or
\[
\frac{h^2}{2\sqrt{h^2-R^2_F}} = \frac{1^2}{2\sqrt{1^2-R^2_F}}\,,
\]
which has two solutions. The first is $h^2 = 1$, yielding (as mentioned already) the
regular $d$-simplex. The second is given by
\[
h = \frac{R_F}{\sqrt{1-R^2_F}}\,.
\]
In order that this value of $h$ gives rise to an honest isosceles
simplex, it is necessary and sufficient that $h > R_T$, that is
\[
\frac{R_F}{\sqrt{1-R^2_F}} > \frac{1}{2\sqrt{1-R^2_F}}
\]
or
\[
R_F > \frac{1}{2}\,.
\]
Since $R_F$, being the circumradius of the regular $(d-2)$-simplex of edge-lenght $1$, increases
with $d$, and since $R_1 = \frac{1}{2}$, it follows that $R_F > \frac{1}{2}$ iff $d-2 > 1$, i.e.,
iff $d \ge 4$. \hfill
$\Box$

\textbf{Remark 2:} One can have an alternative and intuitive view of
the latter construction of equiradial $d$-simplices that are not
equiareal for $d \ge 4$. Namely, consider the regular
$(d-1)$-simplex $T$ of edge length $1$ inscribed in its
circumsphere of radius $R_T$. Now each facet $F$ of $T$ is also
the facet of a second isosceles $(d-1)$-simplex inscribed in the
same sphere, hence having the same circumradius as $T$. These
``ears'' can be folded up to form the desired $d$-simplex provided
the length of the external edges is large enough or, equivalently,
that the height of the ``ears'' is large enough. One can easily
check that $d \ge 4$ suffices to complete the construction. \hfill
$\Box$

We continue with a characterization of regular simplices by two coinciding centers, one of which still
has to be defined. Namely, the \emph{$1$-center} of $S = [A_1, \dots, A_{d+1}]$ is the center of
the $(d-1)$-sphere which is tangent to all edges $A_i A_j$ of $S$ if it exists (in general
it does not exist).

\textbf{Theorem 3.5:} \emph{If the $1$-center of a $d$-simplex $S$ exists and coincides
with the circumcenter of $S$, then $S$ is regular}.

\textbf{Proof:} The circumcenter of $S$ can be viewed as the intersection of the hyperplanes
perpendicular to the edges at their midpoints. On the other hand, the $1$-center is the
intersection of hyperplanes perpendicular to the edges $A_i A_j$ at points dividing any
such edge into lengths $\overline{a}_i$ and $\overline{a}_j$ such that $|A_i A_j| =
\overline{a}_i + \overline{a}_j$, depending only on its endpoints. (Note that all
tangential segments from an exterior point of a $(d-1)$-sphere to the respective touching
points have equal lengths.) The hyperplanes defining the $1$-center are obviously parallel
to the hyperplanes defining the circumcenter. Thus, if the $1$-center exists and coincides with
the circumcenter, the two families of hyperplanes must coincide, and it follows that $\overline{a}_i
= \overline{a}_j = \frac{1}{2} |A_i A_j|$ for all different $i,j \in \{1, \dots, d+1\}$. Hence
all edge lengths of $S$ are equal to $2\overline{a}_i$, i.e., $S$ is regular.
\hfill $\Box$

Finally we mention three characterizations of regular simplices within the restricted family
of orthocentric simplices (since these statements are related to our considerations).
The first one was proved in \cite{PW}: \emph{An orthocentric $d$-simplex
is regular iff its orthocenter and its Fermat-Torricelli point coincide}. In \cite{Fr} it was
shown that \emph{an orthocentric $d$-simplex is regular iff its centroid and its
orthocenter coincide}, and \cite{Ge} contains the observation that \emph{an equiareal orthocentric $d$-simplex
is regular}.

\section{Constructions in dimension \boldmath$4$\unboldmath}

In view of further results restricted to dimension $4$,
we continue with the representation of a tool that relates the
geometry of a simplex $S$ to the algebraic properties of a
certain matrix associated to $S$ (see \cite{HJ} and \cite{LT}).
Namely, for a $d$-simplex $S =
[A_1, \dots, A_{d+1}]$ in $\bE^d$ one defines the Gram matrix
$G(S)$ to be the symmetric, positive semidefinite $(d+1) \times
(d+1)$ matrix of rank $d$ whose $(i,j)$th entry is the inner
product $A_i \cdot A_j$ (we mean the ordinary inner
product, say), cf. \cite{HJ}, p. 407. Given $G(S)$, one can
calculate the distances $d(A_i, A_j)$ for every $i,j$ using the
formula
\[
\left(d(A_i, A_j)\right)^2 = \left(A_i - A_j\right) \cdot
\left(A_i - A_j\right)\,.
\]
According to the last part of Proposition 9.7.1 in \cite{Be},
$G(S)$ determines $S$ up to an isometry of $\bE^d$. Also one
recovers $S$ from $G(S)$ via the Cholesky factorization $G(S) =
HH^t$, where the rows of $H$ are the vectors $A_i$ coordinatized
with respect to some orthonormal basis of $\bE^d$. In fact, if
$G(S)$ is a symmetric, semidefinite, real matrix of rank $r$, say,
then there exists a unique symmetric, positive semidefinite, real
matrix of rank $r$ with $H^2 = G(S)$, cf. \cite{HJ}, Theorem
7.2.6, p. 405, and the symmetry of $H$ implies $G(S) = HH^t$.

A $4$-simplex $S = [A,B,C,D,E]$ whose circumcenter and
centroid coincide (with circumradius $1$ and the origin as circumcenter, say) is
thus determined by unit vectors $A,B,C,D,E$ satisfying $A+B + C +
D + E = o$. The Gram matrix $G(S)$ is then a symmetric, positive
semidefinite matrix of rank 4 having the form
\begin{equation}
G(S) = \left(\begin{array}{lllll} 1 & x & y & z & \bullet\\
\bullet & 1 & Z & Y & \bullet\\
\bullet & \bullet & 1 & X & \bullet\\
\bullet & \bullet & \bullet & 1 & \bullet\\
\bullet & \bullet & \bullet & \bullet & 1 \end{array}\right)\,,
\end{equation}
where $x = A \cdot B, \, y = A \cdot C, \, z = A \cdot D, \, X = C \cdot
D, \, Y = B \cdot D, \, Z = B \cdot C$, and the $\bullet$'s are defined
by the symmetry of $G(S)$ and the fact that the entries of every row add up to
zero. Thus, to construct a simplex $S = [A,B,C,D,E]$ whose
circumcenter and centroid coincide (at $o$, say), we need to
construct a matrix $G(S)$ of the form described in (7) and
satisfying the conditions formulated after (7). We then define
$A,B,C,D, E$ to be the rows of the matrix $H$ that satisfies $HH^t
= G$. To the assumption that the resulting simplex $S$ have $o$ as
its incenter we add the extra requirement that all facets of $S$
have the same $3$-volume. Denoting the $3$-volume of the facet
$[A,B,C,D]$ by $V_E$, we have $4V_E = \det M_E M_E^t$, where $M_E$
is the matrix whose rows are the vectors $B-A, C-A, D-A$. In terms
of $x,y,z, X, Y, Z$ this is written as
\[
\begin{array}{lll}
4V^2_E & = & \det \left(\begin{array}{ccc}
(B-A) \cdot (B-A) & (B-A) \cdot (C-A) & (B-A) \cdot (D-A)\\
\bullet & (C-A) \cdot (C-A) & (C-A) \cdot (D-A)\\
\bullet & \bullet & (D-A) \cdot (D-A)
\end{array}\right)\\
&&\\
& = & \det \left(\begin{array}{ccc}
2-2x & 1-x-y+Z & 1-x-z+Y\\
\bullet & 2-2y & 1-y-z+X\\
\bullet & \bullet & 2-2z \end{array}\right)\,,
\end{array}
\]
where the $\bullet$'s are to be filled in by symmetry. Defining $V_A, V_B, V_C, V_D$ analogously and
recalling that $E = -A-B-C-D$, we get
\[
M_D = \left(\begin{array}{c} B-A\\C-A\\E-A\end{array}\right) =
\left(\begin{array}{c} B-A\\C-A\\-2A-B-C-D\end{array}\right)
\]
\begin{equation}
4V^2_D = \det \left(\begin{array}{ccc}
2-2x & 1-x-y+Z & 1-x+y+z-Z-Y\\
\bullet & 2-2y & 1-y+x+z-X-Z\\
\bullet & \bullet & 4+2x+2y+2z \end{array}\right)\,.
\end{equation}
The matrices $M_C, M_B, M_A$ (and $4V^2_C, 4V^2_B, 4V^2_A$) are obtained from $M_D$ (and
$4V^2_D$) by applying the permutations $(z\,y) (Z\,Y) (z\,x) (Z\,X) (y\,X) (Y\,x)$, respectively.
Thus a simplex $S = [A,B,C,D,E]$ corresponds to a symmetric, positive definite matrix $G(S)$ as in (7)
each of whose rows adds up to zero with
\begin{equation}
4V^2_E = 4V^2_D = 4V^2_C = 4V^2_B = 4V^2_A\,.
\end{equation}

\textbf{Theorem 4.1:} \emph{Let $G$ be a symmetric matrix of the form $(7)$ each row
of which adds up to zero, and let $G_0$ be obtained from $G$ by taking}
\[
y = Y = x\,, \,\, z = Z = X = - \frac{1}{2}-x\,,
\]
\emph{where $x$ is such that
\begin{equation}
\frac{1}{4} < x < \frac{\sqrt{5}-1}{4} \quad
\mbox{$($or, equivalently}, \, 72° < \cos^{-1} x < 144°)\,.
\end{equation}
Let $A,B,C,D,E$ be the row vectors of the matrix $H$ defined by $HH^t = G_0$, and let
$S = [A,B,C,D,E]$. Then the centroid, the circumcenter, the incenter
and the Fermat-Torricelli point of $S$ coincide}.

\textbf{Proof:} The characteristic polynomial of $G_0$
is $T(T-r_1)^2 (T-r_2)^2$ with
\[
r_1, r_2 = \frac{(5 \pm \sqrt{5}) (4x+1)}{4}\,,
\]
as one can immediately check.
By (10), $r_1$ and $r_2$ are positive, and therefore $G_0$ is positive semidefinite
and of rank 4. Thus there exists a real matrix $H$ of rank 4 with $HH^t = G_0$. The
row vectors $A,B,C,D,E$ of $H$ are unit vectors since the diagonal entries of $G_0$ are
$1$'s. From the discussion above we still need only to check that $V_A = V_B = V_C = V_D = V_E$,
which is immediate, with $V_A = \frac{5}{2} \sqrt{1-2x-4x^2}$. This completes the proof.
\hfill $\Box$

\textbf{Corollary 4.2:} \emph{There exists a non-regular $4$-simplex whose centroid,
circumcenter, incenter and Fermat-Torricelli point coincide}.

\textbf{Proof:} This corollary follows from the fact that the simplex $S$ exhibited in Theroem
4.1 is regular iff $x = - \frac{1}{4}$. \hfill $\Box$

\textbf{Theorem 4.3:} \emph{There exists an equiareal $4$-simplex whose centroid, circumcenter
and Fermat-Torricelli point are pairwise distinct}.

\textbf{Proof:} Let
\[
G = \left(\begin{array}{ccccc}
1 & x & x & -1-2x & x\\
x & 1 & 5x & x & x\\
x & 5x & 1 & x & x\\
-1-2x & x & x & 1 & x\\
x & x & x & x & 1 \end{array}\right)\,.
\]
It is routine to check that the characteristic polynomial of $G$ is $g(T) = (T - (2x +2))
(T-(1-5x)) f(T)$, where
\[
f(T) = T^3 - (2 + 3x) T^2 + (1+x-18x^2) T - 2x(x^2 - 8x-1)\,,
\]
and that exactly one of the zeros of $g(T)$ represents $0$ while the others are non-negative
iff $x= 4-\sqrt{17}$. Let $S = [A,B,C,D,E]$ be the $4$-simplex that corresponds to $G$
for this value of $x$. Thus, again $A,B,C,D,E$ are the rows of the matrix $H$ with $G
= HH^t$. Since the diagonal of $G$ consists of $1$'s, the circumcenter of $S$ is ${\bf 0}$.
And since the rows of $G$ do not add to zero, the centroid of $S$ is not ${\bf 0}$. From this and Theorem 4.1
it follows that the centroid, the circumcenter and the Fermat-Torricelli point are pairwise distinct. It remains
to check equiareality. Again it is routine to show that the volumes of all facets of $S$
are equal to $4 - 20x - 4 x^2 + 20x^3$. \hfill $\Box$

On the other hand one might ask how the properties (i) -- (iv) in Theorem 3.2 are connected with each
other. The following statements refer to this in $4$-space.

\textbf{Theorem 4.4:} \emph{For the centroid $\cG$, the circumcenter $\cC$, and the
incenter $\cI$ of a $4$-simplex, the property $\{\cG = \cC\}$
does not imply any of $\{\cC = \cI, \cI = \cG\}$, and the property $\{\cI=\cG\}$ does
not imply any of $\{\cG=\cC, \cC=\cI\}$}.

\textbf{Proof:} The existence of a non-equiareal $4$-simplex with coinciding centroid and
circumcenter was verified in \cite{Ha-Wa 2}. Thus $\cG= \cC \nRightarrow \cC = \cI$ and
$\cG = \cC \nRightarrow \cI = \cG$.
Theorem 4.3 above shows that $\cI = \cG \nRightarrow \cG = \cC$ and $\cI = \cG \nRightarrow \cC = \cI$.
\hfill $\Box$

V. Devide \cite{De 2} asked \emph{whether there are $4$-simplices which are both equiradial and
equiareal, but not equifacetal}. In the following we will answer that question in the
affirmative, by constructing a $1$-parameter family of non-equifacetal $4$-simplices that
are both equiradial and equiareal. Namely, we will consider isosceles $4$-simplices (for this notion
see the proof of Theorem 3.3) whose bases are
equifacetal $3$-simplices $T = [A_1, A_2, A_3, A_4]$ with three distinct edge lengths
$a,b,c$, say. We denote such a tetrahedron by $T = (a,b,c,a,b,c)$, where the edges are written in the
order $a_{12}, a_{23}, a_{13}, a_{34}, a_{14}, a_{24}$, with $a_{ij} = |A_i A_j|$. It is
well known that opposite edges of $T$ are congruent, and that equifacetal tetrahedra exist
iff $a,b,c$ are the edge lengths of an acute triangle, see, e.g., \cite{AC}. This is
characterized by the conditions $a^2 < b^2 + c^2, \, b^2 < a^2 + c^2$, and $c^2 < a^2 + b^2$.
Using the notation from the proof of Theorem 3.3, $T$ is the base of certain isosceles
$4$-simplices $S = [T,P]$, where $P$ lies on a line perpendicular to the affine hull of $T$ at
$T$'s circumcenter. Thus $P$ has to be equidistant to the four vertices of $T$, and so $S =
[T,P]$ can be described, in terms of edge lengths, by $S = (a,b,c,a,b,c,h,h,h,h)$ for an
appropriate $h$. We will show how to choose $a,b,c,h$ for getting $4$-simplices that are equiradial and
equiareal, but not equifacetal. We also mention that these simplices cannot be obtained as
perturbations of regular simplices.

\textbf{Lemma 4.5:} \emph{Let $D$ be an acute $($or right$)$ triangle with side lengths $a,b,c$
and with circumradius $R$. Then}
\[
\frac{a^2 + b^2 + c^2}{8} \ge R^2 \ge \frac{a^2 + b^2 + c^2}{9}\,,
\]
\emph{with the extreme values attained when $D$ is right-angled and when $D$ is equilateral}.

\textbf{Proof:} Let $o$ be the circumcenter of $D$. Applying the Law of Cosines to the triangles $BOC,
COA$, and $AOB$, and using the facts that $\angle BOC = 2A$, etc., we obtain
\[
a^2 + b^2 + c^2 = 6R^2 - 2R^2 \sigma\,,
\]
where $\sigma = \cos 2A + \cos 2B + \cos 2 C$. We have
$\sigma = - 4 \cos A \cos B \cos C - 1$ by \cite{Ca}, formula 682, page 166, and therefore
\begin{equation}
a^2 + b^2 + c^2 = 8R^2 (1 + \cos A \cos B \cos C)\,.
\end{equation}
Since $D$ is acute, it follows that the minimum of $\cos A \cos B \cos C$ is $0$, and is attained
when $D$ is right-angled. Also, the maximum of $\cos A \cos B \cos C$ is $1/8$, and is attained
when $A = B = C$. This follows from the fact that if $x > y$, then
\[
2 \cos x \cos y = \cos (x+y) + \cos (x+y) < 1 + \cos (x+y) = 2 \cos^2 \frac{x+y}{2}\,.
\]
Thus $0 \le \cos A \, \cos B \, \cos C \le 1/8$, with the extreme values attained at right-angled and
equilateral triangles. The rest follows from (11). \hfill $\Box$

\textbf{Theorem 4.6:} \emph{Let $D = (a,b,c)$ be an acute triangle with side lengths $a,b,c$ and
with circumradius $R$. Let $T = (a,b,c,a,b,c)$ be the equifacetal tetrahedron having $D$ as a facet,
and let $S = (a,b,c,a,b,c,h,h,h)$ be the isosceles $4$-simplex obtained by adjoining to
$T$ a vertex at feasible distance $h$ from each vertex of $T$. Then $S$ is equiareal and equiradial iff $S$ is
regular or}
\begin{equation}
R^2 = \frac{3(a^2 + b^2 + c^2)}{25} \, \mbox{ and } \, h^2 = \frac{a^2 + b^2 + c^2}{5}\,.
\end{equation}
\emph{Consequently, there exist equiareal, equiradial $4$-simplices that are not equifacetal}.

\textbf{Proof:} Let $K$ be the area of the triangle $D$ and let $Q = 16K^2$. It is well-known that
\[
\begin{array}{lll}
Q & = & 16K^2 = 2(a^2 b^2 + b^2 c^2 + c^2 a^2)-(a^4 + b^4 + c^4)\\
&&\\
R^2 & = & {\D \frac{a^2 b^2 c^2}{(a+b+c) (-a+b+c) ()a-b+c) (a+b-c)}}\,;
\end{array}
\]
see \cite{Isa}, page 69. Letting
\[
u = a^2 + b^2 + c^2\,, \,\, v = a^2 b^2 + b^2 c^2 + c^2 a^2\,,\,\, w = a^2 b^2 c^2\,,
\]
we see that
\begin{equation}
Q = 4v - u^2 \,, \, \,\, w = R^2 Q\,.
\end{equation}
We find it more convenient to work with the parameters $u, R$, and $Q$ instead of $a,b$, and $c$, and we
freely use the relations in (13).
Recalling that $T$ is the equifacetal tetrahedron $(a,b,c,a,b,c)$, we let $T'$ be the isosceles tetrahedron
$(a,b,c,h,h,h)$ and note that each facet of $S$ other than $T$ is congruent to $T'$.

By the well-known volume formula for $d$-simplices in terms of their edge lengths (see, e.g.,
\cite{PT}, Problem 1.18, page 29), the volume $V$ of the tetrahedron with edge lengths
$x,y,z, X, Y, Z$ is given by
\[
288V^2 = \left[\begin{array}{ccccc}
0 & 1 & 1 & 1 & 1 \\
1 & 0 & x^2 & y^2 & Z^2\\
1 & x^2 & 0 & z^2 & Y^2\\
1 & y^2 & z^2 & 0 & X^2\\
1 & Z^2 & y^2 & X^2 & 0 \end{array}\right]\,.
\]
In particular, with $x = X = a, y = Y =b, z = Z = c$ the volume $V_T$ of $T$ is determined by
\[
\begin{array}{lll}
288V^2_T & = & 4 (b^2 + c^2 - a^2) (c^2 + a^2 - b^2) (a^2 + b^2 - c^2)\\
&&\\
& = & 4 (u-2a^2) (u-2b^2) (u-2c^2) = 4 (-u^3 + 4uv + 8w)\\
&&\\
& = & 4(u(4v-u^2)-8w) = 4(uQ - 8QR^2) = 4Q (u-8R^2)\,,
\end{array}
\]
see \cite{AC}, page 102. Similarly, with $x = a, y= b,
z = c, X = Y = Z = H$ the volume $V_{T'}$ of $T'$ is determined by
\[
\begin{array}{lll}
288V^2_{T'} & = & 4a^2 b^2 h^2 + c^2 b^2 h^2 - 2h^2 b^4 + 4c^2 a^2 h^2 - 2c^4 h^2 - 2a^4 h^2 - 2a^2b^2c^2\\
&&\\
& = & 2h^2 Q - 2w = 2h^2 Q - 2QR^2 = 2Q (h^2-R^2)\,.
\end{array}
\]
Thus equiareality of $S$, given by $V_T = V_{T'}$, is equivalent to the condition
\begin{equation}
h^2 = 2u-15R^2\,.
\end{equation}
Equiradiality of $S$ is equivalent to the condition $R_T = R_{T'}$, where $R_T, R_{T'}$ are the
circumradii of $T$ and $T'$, respectively. By \cite{AC}, p. 102, $R_T$ is given by
\begin{equation}
R^2_T = \frac{a^2 + b^2 + c^2}{8} = \frac{u}{8}\,.
\end{equation}
For $R_{T'}$, we use the formula obtained in the proof of Theorem 3.3, yielding
\[
R^2_{T'} = \frac{h^4}{4(h^2-R^2)}\,.
\]
From this and (15) it follows that equiradiality is equivalent to the condition
\begin{equation}
2h^4 = u (h^2 - R^2)\,.
\end{equation}
It is easy to verify that (12) satisfies both (14) an (16). Conversely, if (14) and (16) hold
then, eliminating $h$, we obtain
\[
2 (2u-15R^2)^2 = u(2u-15R^2-R^2)\,,
\]
which simplifies into
\[
0 = 3u^2 - 52 u R^2 + 225 R^4 = (u-9R^2) (3u - 25R^2)\,.
\]
By Lemma 4.5, the solution $u = 9R^2$ corresponds to the equilateral triangle $a = b = c$.
In view of (14) this corresponds to the case when $h = a$, i.e., to the case when $S$ is the regular
$4$-simplex. The solution $3u = 25R^2$ corresponds, again in view of (14), to the case $h^2 = u/5$.

To prove the last statement, note that if $D = (a,b,c)$ runs over all acute triangles inscribed
in a circle of radius $R$, then, by Lemma 4.5 and continuity arguments, $a^2 + b^2 + c^2$ will take all
values between $8R^2 $ and $9R^2$. Thus, given any $R > 0$, there exists an acute triangle whose side lengths
$a,b, c$ satisfy $25R^2 = 3 (a^2 + b^2 + c^2)$. Finally, to guarantee the existence of the isosceles
$4$-simplex $(a,b,c,a,b,c,h,h,h,h),h$ can take any value that is greater than the circumradius $R_T$
of $T$. Thus, in view of (15), the only restriction on $h$ is given by
\[
h^2 > \frac{a^2 + b^2 + c^2}{8}\,,
\]
and the choice $h^2 = (a^2 + b^2 + c^2)/5$ falls within this restriction.

\textbf{Remark 3:} Note that the construction given above produces $4$-simplices  with $\cC=\cI=\cG \,
(= \cF)$ that are also not equifacetal and hence not regular.

\textbf{Remark 4:} Within the family of isosceles $4$-simplices with an equifacetal base, the degree
of freedom in constructing an equiareal, equiradial, but non-equifacetal simplex is embodied in our
freedom in choosing an acute triangle whose side lengths $a,b,c$ and circumradius $R$ satisfy the
relation
\[
R^2 = \frac{3(a^2 + b^2 + c^2)}{25}\,.
\]
It would be interesting to investigate whether this freedom can be exploited in constructing
$4$-simplices that have, beside equiareality, equiradiality and non-equifacetality, additional significant
properties.

\section{Center coincidence and cevians}

This section will refer to the relationship between coinciding centers on the one hand,
and the lengths of cevians associated with these centers on the other hand.

As already mentioned, the affine independence of the vertex set of $S = [A_1, \dots,
A_{d+1}]$ allows a unique (up to a constant) linear combination of
the origin ${\bf 0}$, i.e.,
\begin{equation}
{\bf 0} = a_1 A_1 + \dots + a_{d+1} A_{d+1}
\end{equation}
with
\begin{equation}
s = a_1 + \dots + a_{d+1} \not= 0\,.
\end{equation}
Namely, otherwise we would have
\[
{\bf 0} = a_1 \left(A_1 - A_{d+1}\right) + \dots + a_d \left(A_d -
A_{d+1}\right)
\]
with $a_1 = \dots = a_d = 0$ and hence $a_{d+1} = 0$,
contradicting the non-triviality of (17). We shall also assume
that none of the vertices of $S$ is ${\bf 0}$, and that the
lines through the vertices and ${\bf 0}$ intersect the
opposite facets. To say that the line joining $A_{d+1}$ and ${\bf 0}$
intersects the opposite facet is equivalent to the existence of
numbers $c_1, \dots , c_{d+1}$ such that
\[
c_1 A_1 + \dots + c_d A_d = c_{d+1} A_{d+1} \, \mbox{ and } \, c_1
+ \dots + c_d = 1\,.
\]
>From the uniqueness of (17) it follows that $a_1 + \dots + a_d
\not= 0$. Therefore we may also assume that
\begin{equation}
\mbox{ no } d \mbox{ of the numbers } \, a_1, \dots, a_{d+1} \, \mbox{ add
up to } \, 0\,.
\end{equation}
Under the assumptions (17), (18), and (19) we let $A^*_j$ be
the point where the line through $A_j$ and ${\bf 0}$ intersects the $j$th
facet of $S$. The line segment $A_j A^*_j$ is usually called the \emph{cevian} through
$A_j$ relative to ${\bf 0}$.
Since
\[
\frac{-a_{d+1} A_{d+1}}{a_1 + \dots + a_d} = \frac{a_1 A_1 + \dots + a_d
A_d}{a_1 + \dots + a_d}\,,
\]
and since the left hand side lies on the cevian and the right
hand side lies on the facet, it follows that
\[
A^*_{d+1} = - \frac{-a_{d+1} A_{d+1}}{a_1 + \dots + a_d}
\]
and
\[
\| A_{d+1} - A^*_{d+1} \| = \frac{|a_1 + \dots + a_{d+1}|}{|a_1 +
\dots + a_d|} \| A_{d+1} \|\,.
\]
Thus we have that
\begin{equation}
\mbox{ the } d+1 \mbox{ cevians through } {\bf 0} \mbox{ are equal }
\Leftrightarrow \frac{|s|}{|s-a_j|} \|A_j\| \, \mbox{ is independent of } \, j\,,
\end{equation}
where $s = a_1 + \dots + a_{d+1}$.

\textbf{Theorem 5.1:} \emph{Let $S = [A_1, \dots, A_{d+1}]$ be a
$d$-simplex. Then the following properties of $S$ are equivalent}.
\begin{enumerate}
\item \emph{The centroid $\cG$ and the circumcenter $\cC$ of $S$ coincide}.
\item \emph{The cevians through the centroid $\cG$ have equal lengths}.
\item \emph{The cevians through the Fermat-Torricelli point $\cF$ have equal lengths}.
\item \emph{The circumcenter $\cC$ lies in $S$ and the cevians through $\cC$
have equal lengths}.
\end{enumerate}

\textbf{Proof:} Assume that the centroid is ${\bf 0}$. Then $a_j = \frac{1}{d+1}$ for all
$j$ in (17). Using (20), we get the following equivalences.
\begin{itemize}
\item[~] ~~ The cevians through the centroid of $S$ are of equal lengths
\item[~] $\Leftrightarrow$ $\frac{d+1}{d} ||A_j ||$ is independent of $j$
\item[~] $\Leftrightarrow$ $||A_1|| = \dots = || A_{d+1} ||$
\item[~] $\Leftrightarrow$ ${\bf 0}$ is the circumcenter
\item[~] $\Leftrightarrow$ the circumenter coincides with the centroid.
\end{itemize}
Assume that the Fermat-Torricelli point is ${\bf 0}$. Then we may suppose that
$a_j = \frac{1}{||A_j||}$ for all
$j$ in (17). Again using (20), we see the following equivalences.
\begin{itemize}
\item[~] ~~ The cevians through the Fermat-Torricelli point of $S$ have equal lengths
\item[~] $\Leftrightarrow$ $a_i (s-a_i) = a_j (s - a_j) \Leftrightarrow (a_i - a_j) (s-a_i-a_j) = 0
\Leftrightarrow a_i = a_j$ for all $i,j$
\item[~] $\Leftrightarrow$ $a_1 = \dots = a_{d+1}$
\item[~] $\Leftrightarrow$ ${\bf 0}$ is the centroid of $S$
\item[~] $\Leftrightarrow$ the Fermat-Torricelli point of $S$ coincides with the centroid of $S$
\item[~] $\Leftrightarrow$ the centroid of $S$ coincides with the circumcenter of $S$ (by Theorem 4.1).
\end{itemize}
Finally, assume that the circumcenter is ${\bf 0}$ and that it lies in $S$. Then in (17) $||A_1|| =
\dots = ||A_{d+1}||$, and $a_j \ge 0$ for all $j$. Using (20), we get the following
equivalences.
\begin{itemize}
\item[~] ~~ The cevians through the circumcenter of $S$ are of equal lengths
\item[~] $\Leftrightarrow$ $|s-a_j|$ is independent of $j$
\item[~] $\Leftrightarrow$ $s-a_j$ is independent of $j$ (because $s - a_j \ge 0$)
\item[~] $\Leftrightarrow $ $a_1 = \dots = a_{d+1}$
\item[~] $\Leftrightarrow$ ${\bf 0}$ is the centroid of $S$
\item[~] $\Leftrightarrow$ the centroid of $S$ and the circumcenter of $S$ coincide. \hfill $\Box$
\end{itemize}
The last statement in Theorem 5.1 triggers the question whether there exists a $d$-simplex
$S$ whose circumcenter $\cC$ is outside of $S$ and whose cevians through $\cC$ have equal
lengths. Theorem 5.5 will show that this may happen if and only if $d \ge 4$. As
preparation for that theorem, we need three lemmas.

\textbf{Lemma 5.2:} \emph{Let $S = [A_1, \dots, A_{d+1}]$ be a $d$-simplex whose circumcenter
is ${\bf 0}$ and whose circumradius is $1$. Then the cevians through the circumcenter are equal if
and only if there exists a suitable $r$ with $0 \le r < \frac{d+1}{2}$ such that, after some
rearrangement, we have}
\begin{equation}
(2d-2r+1) (A_1 + \cdots + A_r) - (2r-1) (A_{r+1} + \dots + A_{d+1}) = {\bf 0}\,.
\end{equation}

\textbf{Proof:} Let the dependence relation among the $A_i$'s be given as in (17), and
let $s = a_1 + \dots + a_{d+1}$. By rearranging the $A_j$'s and by multiplying (17) with $-1$, if
necessary, we may assume that $s-a_j \ge 0$ for $1 \le j \le r$ and $s - a_j < 0$ for $j > r$,
where $0 \le r \le \frac{d+1}{2}$. By (20) we obtain the following equivalences.
\begin{itemize}
\item[~] ~~ The cevians through the circumcenter of $S$ are of equal length
\item[~] $\Leftrightarrow$ $|s-a_j|$ is independent of $j$
\item[~] $\Leftrightarrow$ $s-a_1 = \dots = s-a_r = -(s-a_{r+1}) = \dots = -(s-a_{d+1})$
\item[~] $\Leftrightarrow$ $a_1 = \dots = a_r, \, a_{r+1} = \dots = a_{d+1} =
2s-a_1$
\item[~] $\Leftrightarrow$ $sa_1 = \dots = sa_r = s{\D \frac{2d -2r+1}{d+1-2r}}\,,\, sa_{r+1} =
\dots = sa_{d+1} = s{\D \frac{1-2r}{d+1-2r}}$
\item[~] $\Leftrightarrow$ $(2d-2r+1) (A_1 + \dots +A_r) - (2r-1) (A_{r+1} + \dots +
A_{d+1})={\bf 0}$, as desired.
\end{itemize}
Note that the possibility $d+1 - 2r = 0$ is excluded since it leads to the contradiction
$s = 0$. \hfill $\Box$

\textbf{Lemma 5.3:} \emph{Let $V$ be a unit vector in $\bE^d$, and let $t \not= 0$ be in the
open interval $(-d,d)$. If $d \ge 2$, then there exists a basis $B_1, \dots, B_d$ of~
$\bE^d$ consisting of unit vectors with $B_1 + \dots + B_d = tV$. In fact, one can choose
$B_1, \dots, B_d$ to be equally inclined in the sense that $B_i \cdot B_j = B_k
\cdot B_l$ whenever $i \not= j$ and $k \not= l$}.

\textbf{Proof:} Let $U$ be the orthogonal complement of $V$, and let $[E_1, \dots, E_d]$
be a regular $(d-1)$-simplex in $U$ centred at ${\bf 0}$ and having circumradius $1$. Then
$E_1, \dots, E_d$ are equally inclined, affinely independent unit vectors with $E_1 + \dots +
E_d = {\bf 0}$.
As $x$ takes all non-zero values, $\frac{dx}{\sqrt{1+x^2}}$ takes all non-zero values
in the open interval $(-d, d)$. It follows that there exists an $x$ such that
$\frac{dx}{\sqrt{1+x^2}} = t$. Let
\[
B_j = \frac{E_j + x V}{\sqrt{1+x^2}}\,.
\]
Then the $B_j$'s are equally inclined unit vectors with $B_1 + \dots + B_d = tV$. It remains
to show that they are linearly independent. We use the fact that if a linear combination
of equally inclined unit vectors vanishes, then all the coefficients are equal, see also
\cite{Kl}. Thus
\[
\begin{array}{lll}
c_1 B_1 + \dots c_d B_d = {\bf 0} & \Rightarrow & \left(c_1 E_1 + \dots + c_d E_d\right) + x
\left(c_1 + \dots + c_d\right) V = {\bf 0}\\
& \Rightarrow & c_1 E_1 + \dots + c_d E_d = \left(c_1 + \dots + c_d\right) V = {\bf 0}\\
&& (\mbox{since } V \perp U)\\
& \Rightarrow & c_1 = \dots = c_d \,(\mbox{by \cite{Kl}) and } c_1 + \dots + c_d = 0\\
& \Rightarrow & c_1 = \dots = c_d = 0\,, \mbox{ as desired}.
\quad \quad \quad \quad \quad \quad \quad \quad \quad \quad \quad \quad \quad \Box
\end{array}
\]

\textbf{Lemma 5.4:} \emph{Suppose that $2 \le r \le d-1$ and that $b,c$ are non-zero real numbers
such that $br + c (d-r+1) \not= 0$. Then there exist affinely independent unit vectors
$A_1, \dots, A_{d+1}$ in $\bE^d$ such that}
\begin{equation}
b\left(A_1 + \dots + A_r\right) + c \left(A_{r+1} + \dots + A_{d+1}\right) = {\bf 0}\,.
\end{equation}

\textbf{Proof:} Dividing (22) by an appropriate number, one may assume that $b,c$ are non-zero
small numbers, say in $(-1,1)$. We decompose $\bE^d$ into the direct sum of three
mutually orthogonal subspaces $U_{r-1}, U_1, U_{d-r}$ of dimensions $r-1, 1, d-r$, respectively,
and we let $V$ be a unit vector in $U_1$. By the previous lemma, the direct sum $U_1 \oplus U_{r-1}$
has a basis consisting of unit vectors $A_1, \dots, A_r$ with $A_1 + \dots + A_r = -cV$.
Similarly, $U_1 \oplus U_{d-r}$ has a basis consisting of unit vectors $A_{r+1}, \dots, A_{d+1}$
with $A_{r+1} + \dots + A_{d+1} = bV$. Then $A_1, \dots, A_{d+1}$ satisfy (22). Also, since
the sum of the coefficients in (22) is not zero, being nothing but $br + c (d-r+1)$, it follows
that $o$ is in the affine hull of $A_1, \dots, A_{d+1}$. Hence the affine hull of $A_1, \dots,
A_{d+1}$ is their linear span and thus has dimension $r + (d-r+1)-1=d$. Therefore
$A_1, \dots, A_{d+1}$ are affinely independent, as desired. \hfill $\Box$

\textbf{Theorem 5.5:} \emph{There exists a $d$-simplex whose circumcenter is exterior and
whose cevians through the circumcenter are of equal lengths if and only if $d \ge 4$}.

\textbf{Proof:} We use Lemma 5.2. Taking $d = 2$ and $r = 0$ in (21), we obtain $A_1 + A_2 + A_3 = {\bf 0}$, and
hence ${\bf 0}$ is the centroid and cannot be exterior. Taking $d = 2$ and $r =1$, we obtain $3A_1 -
(A_2 + A_3) = {\bf 0}$, which is impossible since $||3A_1||=3$.
Similarly for the cases $d = 3, r = 0$ and $d = 3,
r = 1$. For $d \ge 4$, let $r$ be any number such that $2 \le r < \frac{d+1}{2}$, and let $b = 2d
-2r + 1\,, \, c = -(2r-1)$. Then $r \le d-1$ and $br + c (d-r+1)\not= 0$. Hence, by Lemma 5.4
there exist affinely independent unit vectors $A_1, \dots, A_{d+1}$ such that $b(A_1 + \dots
+ A_r) + c(A_{r+1} + \dots + A_{d+1}) = {\bf 0}$. Now the simplex $S = [A_1, \dots, A_{d+1}]$ has the desired
properties. \hfill $\Box$

Unfortunately, we were not able to prove a statement in the spirit of Theorem 5.1 that
refers to cevians of equal lengths going through the
incenter of $S$. Such a result would be the natural generalization of the
well-known Steiner-Lehmus theorem, see \cite{Co}, pages 9 and 420.


\begin{thebibliography}{00000}
\bibitem{AC}{\sc N. Altshiller-Court:} \emph{Modern Pure Solid Geometry}. Chelsea Publishing
Co., New York, 1964.
%
\bibitem{Ave} {\sc G. Averkov:} \emph{On the geometry of simplices in Minkowski spaces}.
Stud. Univ. Zilina Math. Ser. \textbf{16} (2003), 1-14.
%
\bibitem{Be}{\sc M. Berger:} \emph{Geometry I}. Springer-Verlag, Berlin, 1994.
%
\bibitem{BMS}{\sc V. Boltyanski, H. Martini, V. Soltan:} \emph{Geometric Methods and
Optimization Theory}. Kluwer Academic Publishers, Dordrecht, 1999.
%
\bibitem{Br}{\sc B. H. Brown:} \emph{A theorem of Bang, isosceles tetrahedra}. Amer. Math.
Monthly \textbf{33} (1926), 224-226.
%
\bibitem{Ca}{\sc G. S. Carr:} \emph{Formulas and Theorems in Mathematics}. Chelsea Publishing
Co., 1970.
%
\bibitem{CB}{\sc P. Couderc, A. Ballicioni:} \emph{Premier livre du tètraèdre}. Gauthier-Villars,
Paris, 1953.
%
\bibitem{Co}{\sc H. S. M. Coxeter:} \emph{Introduction to Geometry}. 2nd ed., John Wiley and
Sons, Inc., New York 1969.
%
\bibitem{Da}{\sc L. Dalla:} \emph{A note on the Fermat-Torricelli point of a $d$-simplex}. J.
Geom. \textbf{70} (2001), 38-43.
%
\bibitem{De 1}{\sc V. Devidé:} \emph{Über gewisse Klassen von Simplexen}. Rad Jugoslav.
Akad. Znanost. Umjetnost. \textbf{370} (1975), 21-37.
%
\bibitem{De 2}{\sc V. Devidé:} \emph{Über eine Klasse von Tetraedern}. Rad Jugoslav. Akad.
Znanost. Umjetnost. \textbf{408} (1984), 45-50.
%
\bibitem{Ed}{\sc A. L. Edmonds:} \emph{The geometry of an equifacetal simplex}.
Indiana University reprint 2003, Available as arXiv: math. MG/0408132.
%
\bibitem{F-M}{\sc P. Frankl, H. Maehara:} \emph{Simplices with given $2$-face areas}. European
J. Combin. \textbf{11} (1990), 241-247.
%
\bibitem{Fr}{\sc R. Fritsch:} \emph{Höhenschnittpunkte für $n$-Simplizes}. Elem. Math.
\textbf{31} (1976), 1-8.
%
\bibitem{Ge}{\sc L. Gerber:} \emph{The orthocentric simplex as an extreme simplex}.
Pacific J. Math. \textbf{56} (1975), 97-111.
%
\bibitem{Ha 1}{\sc M. Hajja:} \emph{A vector proof of a theorem of Bang}.
Amer. Math. Monthly \textbf{108} (2001), 562-564.
%
\bibitem{Ha-Wa 1}{\sc M. Hajja, P. Walker:} \emph{Equifacial tetrahedra}. Internat. J. Math.
Ed. Sci. Tech. \textbf{32} (2001), 501-508; Corrigendum: \emph{Equifaciality of tetrahedra whose incenter
and Fermat-Torricelli point coincide}. Preprint, submitted.
%
\bibitem{Ha-Wa 2}{\sc M. Hajja, P. Walker:} \emph{The measure of solid angles in $n$-dimensional
Euclidean space}. Internat. J. Math. Ed. Sci. Tech. \textbf{32} (2002), 725-729.
%
\bibitem{H-W}{\sc H. Havlicek, G. Weiss:} \emph{Altitudes of a tetrahedron and traceless
quadratic forms}. Amer. Math. Monthly \textbf{110} (2003), 679-693.
%
\bibitem{Ho}{\sc R. Honsberger:} \emph{Mathematical Gems II}. Dolciani Mathematical
Expositions No. 2., Mathematical Association of America, Washington, D.C.,  1976.
%
\bibitem{HJ}{\sc R. A. Horn, C. R. Johnson:} \emph{Topics in Matrix Analysis}. Cambridge University Press,
Cambridge, 1994.
%
\bibitem{Isa}{\sc I. M. Isaacs:} \emph{Geometry for College Students}. Brooks/Cole,
U.S.A., 2001.
%
\bibitem{Ki}{\sc C. Kimberling:} \emph{Triangle Centers and Central Triangles}. Congressus
Numerantium, Vol. 29, Utilitas Mathematica Publishing Incorporated, Winnipeg, Canada, 1998.
%
\bibitem{Kl}{\sc M. S. Klamkin:} \emph{On some symmetric sets of unit vectors}. Math. Mag.
\textbf{64} (1991), 271-273.
%
\bibitem{KM 1}{\sc Y. S. Kupitz, H. Martini:} \emph{The Fermat-Torricelli point and isosceles
tretrahedra}. J. Geom. \textbf{49} (1994), 150-162.
%
\bibitem{KM 2}{\sc Y. S. Kupitz, H. Martini:} \emph{Geometric aspects of the generalized
Fermat-Torricelli problem}. In ``Intuitive Geometry'', Eds. I. Bárány and K. Böröczky,
Bolyai Soc. Math. Studies \textbf{6} (1997), 55-127.
%
\bibitem{LT}{\sc P. Lancaster, M. Tismenetsky:} \emph{The Theory of Matrices}. 2nd ed.,
Academic Press, New York, 1985.
%
\bibitem{Le}{\sc H. Lenz:} \emph{Über einen Satz von June Lester zur Charakterisierung
euklidischer Bewegungen}. J. Geom. \textbf{28} (1987), 197-201.
%
\bibitem{Man}{\sc S. R. Mandan:} \emph{Altitudes of a simplex in $n$-spaces}. J. Austral.
Math. Soc. \textbf{2} (1961/62), 403-424.
%
\bibitem{Ma}{\sc H. Martini:} \emph{A hierarchical classification of Euclidean polytopes with regularity
properties}. In: Polytopes -- Abstract, Convex and Computational, Eds. T. Bisztriczky,
P. McMullen, R. Schneider, and A. Ivi\`c-Weiss, Nato ASI Series, Ser. C: Mathematical and
Physical Sciences, Vol. 440, Kluwer Acad. Publ., Dordrecht 1994, 71-96.
%
\bibitem{Ma-Wei}{\sc H. Martini, B. Weissbach:} \emph{Napoleon's theorem with weights
in $n$-space}. Geom. Dedicata \textbf{74} (1999), 213-223.
%
\bibitem{M-W}{\sc H. Martini, W. Wenzel:} \emph{Simplices with congruent $k$-faces}. J. Geom.
\textbf{77} (2003), 136-139.
%
\bibitem{McM}{\sc P. McMullen:} \emph{Simplices with equiareal faces}. Discrete Comput. Geom.
\textbf{24} (2000), 397-411.
%
\bibitem{Mo}{\sc F. Molnár:} \emph{Über die Eulersche Gerade und die Feuerbachsche Kugel
des $n$-dimensionalen simplexes} (Hungarian). Mat. Lapok \textbf{11} (1960), 68-74.
%
\bibitem{PW}{\sc C. M. Petty, D. Waterman:} \emph{An extremal theorem for $n$-simplices}.
Monatsh. Math. \textbf{59} (1955), 320-322.
%
\bibitem{PT}{\sc V. V. Prasolov, V. M. Tikhomirov:} \emph{Geometry}. Translations of
Mathematical Monographs, Vol. 200, Amer. Math. Soc., Providence, R. I., 2001.
%
\bibitem{Th}{\sc V. Thébault:} \emph{Géométrie dans l'espace (Géométrie du tétraèdre)}. Libraire
Vuibert, Paris, 1956.
%
\bibitem{We}{\sc B. Weissbach:} \emph{Euklidische $d$-Simplexe mit inhaltsgleichen $k$-Seiten.}
J. Geom. \textbf{69} (2000), 227-233.
%
\bibitem{Za}{\sc M. Zacharias:} \emph{Elementargeometrie und elementare nicht-euklidische
Geometric in synthetischer Behandlung}. In: Encykl. Math. Wiss III, 1. Teil, 2. Hälfte, Leipzig
1913, \S~21.
\end{thebibliography}
\end{document}